\documentclass[dvips,11pt]{article}
\usepackage{amsmath,amssymb,amsfonts,amscd,amsthm,rotating,pstricks,pst-node,bm,xspace}
\def\simarrow{\mathrel{\raise -0.5mm\hbox{$\sim$}}\hspace{-1.8mm}{\rightarrow} } 

\def\bsimarrow{\leftarrow\hspace{-0.7mm}\mathrel{\raise -0.5mm\hbox{$\backsim$}} }

\def\bt{\begin{tabular}}
\def\te{\end{tabular}}

\def\lettrine#1#2#3{\noindent\hangindent#1\hangafter-#2
\hskip-#1\smash{\hbox to #1{#3\hfill}}\ignorespaces}

\newcommand{\To}[1]{\mathop{\to}\limits_{#1}}

\def\BM{\begin{pmatrix}}
\def\EM{\end{pmatrix}}

\def\d=f{\buildrel\hbox{\scriptsize d\'{e}f}\over \Longleftrightarrow}

\def\rit{\text{\it I\hskip -2pt  R}}

\def\nit{\text{\it I\hskip -2pt  N}}

\def\Bd{{\text B}}

\def\Ed{{\text E}}

\def\be{\begin{equation}}
\def\ee{\end{equation}}
\def\beqn{\begin{eqnarray}}
\def\eeqn{\end{eqnarray}}
\def\nobeqn{\begin{eqnarray*}}
\def\noeeqn{\end{eqnarray*}}
\def\ba{\left(\begin{array}}
\def\ea{\end{array} \right) }

\def\u{\underline}
\def\o{\overline}
\def\and{\; \mbox{and} \;}

\newcommand{\half}{\frac{1}{2}}

\def\hfl#1#2{\smash{\mathop{\hbox to 12mm{\rightarrowfill}}
\limits^{\scriptstyle #1}_{\scriptstyle #2}}}

\def\mod{\mathop{\rm mod}\nolimits}

\def\Be{\begin{enumerate}}
\def\Ee{\end{enumerate}}

\def\Bena{\begin{enumerate}
\def\labelenumi{\theenumi)}
\def\theenumi{\arabic{enumi}}
\def\labelenumii{\theenumii)}
\def\theenumii{\alph{enumii}}}

\def\Bean{\begin{enumerate}
\def\labelenumii{\theenumii)}
\def\theenumii{\arabic{enumii}}
\def\labelenumi{\theenumi)}
\def\theenumi{\alph{enumi}}}

\def\Bero{\begin{enumerate}
\def\labelenumii{\theenumii)}
\def\theenumii{\arabic{enumii}}
\def\labelenumi{(\theenumi)}
\def\theenumi{\roman{enumi}}}

\def\BeRo{\begin{enumerate}
\def\labelenumii{\theenumii)}
\def\theenumii{\arabic{enumii}}
\def\labelenumi{(\theenumi)}
\def\theenumi{\Roman{enumi}}}

\def\Bi{\vskip 11pt\begin{itemize}\itemsep=18pt}
\def\bi{\begin{itemize}}
\def\Ei{\end{itemize}\vskip 11pt}
\def\ei{\end{itemize}}

\def\Bd{\begin{description}}
\def\Ed{\end{description}}

\def\R{\right}
\def\L{\left}

\usepackage[latin1]{inputenc}

\def\bsm{(bisemi)}
\def\ST{structures\xspace}
\def\smst{semistructure\xspace}
\def\SMST{semistructures\xspace}
\def\bsmst{bisemistructure\xspace}
\def\BSMST{bisemistructures\xspace}
\def\bigoplus{\mathop{\oplus}\limits}

\def\prod{\mathop{\Pi}\limits}
\def\sum{\mathop{\Sigma}\limits}

\def\bbf{\boldmath\bf}
\def\o{\overline}

\def\Bi{\begin{itemize}}
\def\Ei{\end{itemize}}

\newcommand{\ZZ}{\mathbb{Z}\,}
\newcommand{\NN}{\mathbb{N}\,}

\newcommand{\Aa}{\mathbb{A}\,}
\newcommand{\CC}{\mathbb{C}\,}

\newcommand{\QQ}{\mathbb{Q}\,}
\newcommand{\FF}{\mathbb{F}\,}

\def\det{\operatorname{det}}
\def\Gal{\operatorname{Gal}}
\def\Frob{\operatorname{Frob}}

\def\mod{\operatorname{mod}}

\def\Hom{\operatorname{Hom}}

\def\cusp{\operatorname{cusp}}
\def\Irr{\operatorname{Irr}}

\def\Rep{\operatorname{Rep}}

\def\Ind{\operatorname{Ind}}

\def\GL{\operatorname{GL}}

\def\sym{\operatorname{sym}}

\def\bM{\begin{matrix}}
\def\eM{\end{matrix}}

\def\lr{left (resp. right) }

\def\resp#1{(resp. #1)}
\def\rresp#1{\qquad \mbox{(resp.} \quad #1\ )}

\def\To{\begin{CD} @>>>\end{CD}}

\def\RL{_{R\times L}}

\newtheorem{propo}{\hspace{-4.5mm}{$\bullet$}\ Proposition}[section]
\newtheorem{coro}[propo]{\hspace{-4.5mm}{$\bullet$}\ Corollary}
\newtheorem{Mpropo}[propo]{\hspace{-4.5mm}{$\bullet$}\ Main Proposition}

\begin{document}

\setcounter{page}{0}
{\pagestyle{empty}
\null\vfill
\begin{center}
{\LARGE The Langlands functoriality conjecture in the bisemialgebra framework}
\vfill
{\sc C. Pierre\/}
\vskip 11pt

Institut de Mathématique pure et appliquée\\
Université de Louvain\\
Chemin du Cyclotron, 2\\
B-1348 Louvain-la-Neuve,  Belgium\\
pierre@math.ucl.ac.be

\vfill
Mathematics subject classification (2000): 11G18, 11R34, 11R37, 11R39.
\vfill

\begin{abstract}
The  Langlands functoriality conjecture envisaged in the bisemistructure framework is proved to correspond to the non-orthogonal completely reducible cuspidal representation of bilinear algebraic semigroups.

\end{abstract}
\vfill
\eject
\end{center}

\vfill\eject}
\setcounter{page}{1}
\def\thepage{\arabic{page}}
{\parindent=0pt 
\include{Functor0}
\section{Historical frame of the Langlands functoriality}

The Langlands program originated from {\bbf Artin's reciprocity law of abelian class fied theory\/}.
The simplest version of Artin's reciprocity law states that if $\sigma :\Gal(E/\QQ)\to\CC^*$ is the homomorphism from the Galois group of a finite extension $E$ of $\QQ$ into $\CC^*$, then there exists a Dirichlet character $\chi _\sigma (\ZZ/N\ \ZZ)^*\to\CC^*$ such that
$\sigma (\Frob p)=\chi _\sigma (p)$ for all primes $p$ (unramified) in $E$.

The non abelian equivalent of this reciprocity law concerns the $n$-dimensional representations of the Galois group $\Gal(E/\QQ)$
throughout homomorphisms:
\[\sigma :\quad \Gal(E/\QQ)\To \GL_n(\CC)\]
which led Artin to introduce $L$-functions:
\[ L(s,\sigma )=\prod_p(\det [I_n-\sigma (\Frob p)\ p^{-s}])^{-1}\]
where $\det[I_n-\sigma (\Frob p)\ p^{-s}]$ are characteristic polynomials related to conjugacy classes $A_p$ of $\Gal(E/\QQ)$ described by $\sigma (\Frob p)$ \cite{Gel}, \cite{Rog}.  But, he did not find the $n$-dimensional analogues of Dirichlet characters and $L$-functions.
\vskip 11pt

It was {\bbf Langlands\/} \cite{Lan1} who {\bbf generalizes the concept of Dirichlet characters by introducing\/}:
\Bena
\item a (unitary) cuspidal automorphic representation $\pi $ of $\GL_n(\Aa_F)$, where $\Aa_F$ is the ring of adeles of the global number field $F$ of characteristic zero.

\item the $L$-function associated with $\pi $.
\Ee
\vskip 11pt

These considerations where then formulated in the {\bbf Langlands (global) reciprocity conjecture\/} which asserts that:

\begin{quote}
For any irreducible representation $\sigma $ of $\Gal(\o F/F)$ in $\GL_n(\CC)$, there exists a cuspidal automorphic representation $\pi $ of $\GL_n(\Aa_F)$ in such a way that the Artin $L$-function of $\sigma $ agrees with the Langlands $L$-function of $\pi $ at almost every place where $\pi $ is unramified \cite{Kna}.
\end{quote}
\vskip 11pt

{\bbf In the local case\/}, i.e. when the envisaged field is  a finite extension $K$ of $\QQ_p$, {\bbf Langlands conjectured the existence of bijections between the set $WD\Rep_n(W_K)$\/} of equivalence classes of $n$-dimensional Frobenius semisimple Weil-Deligne representations of the Weil group $W_K$ {\bbf and the set $\Irr\cusp(\GL_n(K))$\/} of equivalence classes of (super)cuspidal irreducible representations of $\GL_n(K)$.
\vskip 11pt

{\bbf The global correspondences of Langlands over function fields \/} on a smooth curve over $\FF_q$ were extensively studied by L. Lafforgue \cite{Laf} while {\bbf the local correspondences of Langlands over finite number fields\/} were worked out by G. Henniart \cite{Hen} and by M. Harris and R. Taylor \cite{H-T}, \cite{Cara}.
\vskip 11pt

One of the major innovations introduced by Langlands in this context is the {\bbf existence of $L$-groups $^LG$\/} consisting in the semidirect product $\widehat G\rtimes \Gal(\o F/F)$ of the complex reductive group $\widehat G$ by the Galois groups $\Gal(\o F/F)$ \cite{Lan2}.
\vskip 11pt

This allows to tie the cuspidal automorphic representation $\pi $ of a reductive group $G(\Aa_F)$, as well as the associated cuspidal representation $\Pi $ of  $\GL_n(\Aa_F)$, to the $n$-dimensional holomorphic representation $\rho $ of the corresponding $L$-group $^LG$.
\vskip 11pt

According to \cite{Lan3}, if $A_v(\pi )$ denotes the $v$-th cuspidal conjugacy class of $G(\Aa_F)$ and $A_v(\Pi )$ the $v$-th cuspidal conjugacy class of $\GL_n(\Aa_F)$, the following equalities can be stated:
\Bean
\item \quad $\{A_v(\Pi )\}_v=\{\rho (A_v(\pi )\}_v$, $\forall\ v\in\NN$, prime,
\item \quad $L(s,\pi ,\rho )=L(s,\Pi )$,\Ee
relating the $L$-function $L(s,\Pi )$ on $\GL_n(\Aa_F)$ to the $L$-function $L(s,\pi ,\rho )$ on $G(\Aa_F)$ and $^LG$.

This shows \cite{Lau} that {\bbf the $L$-group $^LG$ can be interpreted as a dual group of $G(\Aa_F)$\/} in such a way that:
\Bean
\item the cuspidal representation $\rho $ of $^LG$ be the contragradient representation $\widehat \pi $ with respect to the cuspidal representation $\pi $ of $G(\Aa_F)$ \cite{Cara};
\item the coroots of $\widehat \pi $ correspond to the roots of $\pi $ \cite{J-L}.
\Ee
\vskip 11pt

To proceed further into the program of Langlands, the {\bbf principle of functoriality\/} must be envisaged.  In its general form, it can be stated as follows:
\vskip 11pt

Let $G$ and $G'$ be two reductive groups and let $^LG$ and $^LG'$ be the corresponding $L$-groups.

The two homomorphisms:
\begin{alignat*}{3}
 \phi _G &: \quad &G  &\To   G'\\
\text{and} \qquad  \phi _{^LG}&: \quad &^LG   &\To  ^LG'\end{alignat*}
then imply the corresponding homomorphisms:
\begin{alignat*}{3}
 \pi _{\phi _G} &: \quad &\pi =\otimes\pi _v  &\To   \pi '=\otimes\pi '_v\\
\text{and} \qquad  \pi _{\phi _{^LG}}&: \quad &\widetilde \pi =\otimes\widetilde \pi _v  &\To  \widetilde \pi' =\otimes\widetilde \pi' _v\end{alignat*}
of their respective cuspidal representations in such a way that:
\[ \pi '_v = \phi _{G|v}(\pi _v)\; ,\qquad
 \widetilde  \pi '_v = \phi _{^LG|v}(\widetilde \pi _v)\; .\]
 
More concretely, let $\sym^m:\GL_2(\Aa_F)\to \GL_{m+1}(\Aa_F)$ be the $m$-th symmetric power representation of $\GL_2(\Aa_F)$.  If $\pi =\otimes\pi _v$ is the cuspidal representation of $\GL_2(\Aa_F)$, then:
\[\sym^m(\pi )=\otimes\sym^m(\pi _v)\]
will be the searched cuspidal representation of $\GL_{m+1}(\Aa_F)$ in such a way that $\sym^m(\pi_v )$ be the $v$-th cuspidal conjugacy class of $\GL_{m+1}(F_v)$.  This constitutes the {\bbf Langlands functoriality conjecture\/} extensively studied now.

It was especially proved by H. Kim and F. Shahidi \cite{C-K-P-S}, \cite{K-S} by  $L$-function methods for $m=3,4$ when $F$ is a number field.
\vskip 11pt
In order to try a new breakthrough in this problem of functoriality, {\bbf the \bsmst framework\/} introduced in \cite{Pie1} {\bbf will be considered in this paper\/}.

Indeed, it was shown in \cite{Pie2} that, {\bbf to each mathematical structure, it can be generally associated a triple\/} (left \smst, right \smst, \bsmst) where the left and right \SMST, referring (or localised in) respectively to the upper and lower half spaces, operate on each other by means of their product giving rise to a \bsmst in such a way that its bielements be either diagonal bielements \cite{Lang} or cross bielements.
\vskip 11pt

{\bbf Taking into account \BSMST allows:}
\Bi
\item {\bbf to generate richer mathematical \ST\/} by the consideration of cross \bsm\-sub\ST;
\item {\bbf to have a better visibility of the endomorphisms\/} of these (bisemi)structures of which objects are in some cases de facto bilinear, as for instance (square) matrices $(n\times n)$.
\Ei
\vskip 11pt
Before taking up functoriality in this new bilinear mathematical frame, I shall introduce or recall the structure of the considered algebraic groups and of their cuspidal representations in this new context as well as the connection of these with their classical correspondents.

This will constitute the contents of the next two sections.
\vskip 11pt

\section{Structure of the algebraic bilinear semigroups}

\Bi
\item In the framework of \BSMST, {\bbf a group $G$ must be viewed as a triple\linebreak $(G_L,G_R,G\RL)$ where\/}:
\Bi
\item $G_L$ \resp{$G_R$} is a \lr semigroup under the addition of its \lr elements $g_{L_i}$ \resp{$g_{R_i}$} referring to the upper \resp{lower} half space;
\item {\bbf $G\RL$ is a bisemigroup, or bilinear semigroup\/}, such that its bielements $(g_{R_i}\times g_{L_i})$ be submitted to the cross binary operation $\u\times$ defined by:
\[\begin{matrix}
G\RL & \u\times & G\RL & \quad \To \quad & G\RL\\
(g_{R_i}\times g_{L_i}) & \u\times & 
(g_{R_j}\times g_{L_j})& \quad \To \quad & 
(g_{R_i}+                                                                   g_{R_j})\times (g_{L_i}+ g_{L_J})\;.
\end{matrix}\]
By this way, pairs of right and left elements are sent either in diagonal bielements $(g_{R_i}\times g_{L_i})$ and $(g_{R_j}\times g_{L_j})$ or in cross bielements $(g_{R_i}\times g_{L_j})$ and
$(g_{R_j}\times g_{L_i})$.
\Ei

\item Thus, an algebraic group $\GL_n(F)$ of $(n\times n)$ square invertible matrices with entries in a field $F$ can be viewed as a bilinear algebraic semigroup $\GL_n(F_R\times F_L)$ with entries in the product $F_R\times F_L$ of the right and left finite algebraic symmetric extensions $F_R$ and $F_L$ of a global number field $k$ of characteristic zero.

{\bbf This bilinear algebraic semigroup $\GL_n(F_R\times F_L)$ can then be decomposed\/} into the product of the right semigroup $T_n^t(F_R)$ of lower triangular matrices with entries in the semifield $F_R$ by the left semigroup $T_n(F_L)$ of upper triangular matrices with entries in the semifield $F_L$ {\bbf according to:
\[\GL_n(F_R\times F_L) = T^t_n(F_R) \times T_n(F_L)\;.\]}
The \lr algebraic semigroup $T_n(F_L)$ \resp{$T^t_n(F_R)$} can be viewed as an operator:
\begin{align*}
T_n(F_L) : \quad F_L & \To T^{(n)}(F_L)\equiv V_L\\
\rresp{
T^t_n(F_R) : \quad F_R & \To T^{(n)}(F_R)\equiv V_R},\end{align*}
sending the \lr semifield $F_L$ \resp{$F_R$} into {\bbf the \lr $T_n(F_L)$-semimodule $T^{(n)}(F_L)$ \resp{$T^t_n(F_R)$-semimodule $T^{(n)}(F_R)$}} which {\bbf is a \lr affine semispace $V_L$ \resp{$V_R$}} of dimension $n$ localized in the upper \resp{lower} half space.

\item The left semigroup $T_n(F_L)$ then corresponds to the parabolic subgroup $P_n(F)$ of Borel upper triangular matrices over the number field $F$ in such a way that {\bbf the homomorphism\/}:
\[ \Ind_{T\to G} : \quad T_n(F_L) \To \GL_n(F_R\times F_L)\]
{\bbf produces an induction \cite{Rod} from the parabolic subsemigroup $T_n(F_L)$ to the algebraic bilinear semigroup $\GL_n(F_R\times F_L)$\/}.

\item Therefore, it appears that if $T_n(F_L)$ is assumed to be a semisimple (reductive) (semi)group $G$, then {\bbf the associated (semi)group $T_n^t(F_R)$\/}, on which a contragradient (cuspidal) representation can be defined, is the dual (semi)group of $G$ and {\bf would correspond to the Langlands $L$-group $^LG$ of $G$\/} if $^LG$ is interpreted as the dual group of $G$ as suggested in section 1.

\begin{propo} The representation (bisemi)space of the algebraic bilinear semigroup of matrices $\GL_n(F_R\times F_L)$ is given by the $\GL_n(F_R\times F_L)$-bisemimodule $G^{(n)}(F_R\times F_L)$ which is a $n^2$-dimensional affine bisemispace $(V_R\otimes_{F_R\times F_L}V_L)$ belonging to a neutral Tannakian tensor category ${\u{c}}\RL$ equivalent to the category of finite dimensional representations of affine bisemigroup schemes.
\end{propo}

\begin{proof} \Be
\item It was proved in \cite{Pie1} and in \cite{Pie2} that the $n^2$-dimensional affine bisemispace $(V_R\otimes_{F_R\times F_L}V_L)$ is a $\GL_n(F_R\times F_L)$-bisemimodule $G^{(n)}(F_R\times F_L)$ under the action right by left of $\GL_n(F_R\times F_L)$ on an irreducible (unitary) $n^2$-dimensional affine bisemispace in such a way that $G^{(n)}(F_R\times F_L)=T^{(n)}(F_R)\otimes T^{(n)}(F_L)$ be the tensor product $(V_R\otimes_{F_R\times F_L}V_L)$ of the right and left affine semispaces $V_R$ and $V_L$ over the bisemifield $F_R\times F_L$.

\item A \lr algebraic semigroup $T_n(F_L)$ \resp{$T^t_n(F_R)$} is in fact a \lr semigroup scheme over $k$, i.e. a representable functor from the quotient $k$-algebra $Q_L$ \resp{$Q_R$} of the polynomial ring $k[x]$ modulo the ideal $I_L$ \resp{$I_R$} to the \lr affine semispace $V_L$ \resp{$V_R$} \cite{Pie1}.

And, $\GL_n(F_R\times F_L)$ is an affine bisemigroup scheme, i.e. a representable (bi)functor from $Q_R\otimes Q_L$ to the affine bisemispace $(V_R\otimes_{F_R\times F_L}V_L)$.

Then, the category of finite dimensional representations of affine bisemigroup schemes is equivalent in the perspective of \cite{D-M} to the neutral Tannakian tensor category ${\u{c}}\RL$ for which there exists a $F_R\times F_L$-(bi)linear tensor functor $\omega \RL:{\u{c}}\RL\to {\u{G}}^{(n)}(F_R\times F_L)$ where${\u{G}}^{(n)}(F_R\times F_L)$ denotes the category of $n^2$-dimensional bisemimodules $G^{(n)}(F_R\times F_L)$ over $(F_R\times F_L)$.\qedhere\Ee\end{proof}

\item Considering the Artin's reciprocity law recalled in section 1 and the congruence subgroups used in the theories of cuspidal forms, (pseudo-ramified) complex completions of symmetric finite (closed) algebraic extensions of the number field $k$ of characteristic zero were defined \cite{Pie1} at infinite complex places by their degree given by the integers modulo $N$.  They correspond to $\Hom_k(F_L,\CC)$ and $\Hom_k(F_R,\CC)$ where $k$ may be $\QQ$ \cite{Kna}.

By this way, {\bbf we get a \lr tower}
\begin{align*}
F_\omega  &=\{F_{\omega _1},\dots,F_{\omega _{j,m_j}},\dots,F_{\omega _r},F_{\omega _{r,m_r}}\}\\
\rresp{
F_{\o \omega}  &=\{F_{\o\omega _1},\dots,F_{\o\omega _{j,m_j}},\dots,F_{\o\omega _r},F_{\o\omega _{r,m_r}}\}}\end{align*}
{\bf of packets of equivalent complex completions localized in the upper \resp{lower} half space and characterized by their degrees\/} (or ranks):
\begin{align*}
[F_{\omega _j}:k] &= *+m^{(j)}\cdot j\cdot N\;, && j\in\nit\;,\\
\rresp{
[F_{\o\omega _j}:k] &= *+m^{(j)}\cdot j\cdot N} && 1\le j\le r\le\infty \;,
\end{align*}
at every complex place $\omega _j$ \resp{$\o\omega _j$}

where:
\Bi
\item the closure ${\o{\omega _j(F)}}$ of $\omega _j(F)$ is  a compact field;
\item $*$ denotes an integer inferior to $N$;
\item $m^{(j)}$ is the number (or multiplicity) of the real completions covering $F_{\omega _j}$ \resp{$F_{\o\omega _j}$};
\item $j$ is a global residue degree $f_{\omega _j}$ \resp{$f_{\o\omega _j}$}.
\Ei

When the global Weil groups are considered, the degrees of the finite algebraic extensions $\widetilde F_{\omega _j}$ \resp{$\widetilde F_{\o\omega _j}$} associated with the completions $F_{\omega _j}$ \resp{$F_{\o\omega _j}$} are restricted to \cite{Pie1}:
\begin{align*}
[\widetilde F_{\omega _j}:k] &\equiv \quad [\widetilde F_{\o\omega _j}:k]=0\ \mod N\\
\noalign{\qquad \quad leading to}
[\widetilde F_{\omega _j}:k] &\equiv \quad [\widetilde F_{\o\omega _j}:k]=j\cdot m^{(j)}\cdot N\;.\end{align*}

Each complex \resp{conjugate complex} completion $F_{\omega _p}$ \resp{$F_{\o\omega _p}$} at the primary complex infinite place $\omega _p$ \resp{$\o\omega _p$}, associated with the corresponding finite extension of $k$, is then characterized by a degree: 
\[ [F_{\omega _p};k]\equiv \quad [F_{\o\omega _p}:k]=m^{(p)}\cdot p\cdot N\]
and a number of elements:
\[ n_e (p)=m^{(p)}\cdot p\cdot N\cdot (n_{\rm nu})\]
where $(n_{\rm nu})$ is the number of nonunits.

If the global residue degree is given by
\[ f_{\omega _{p+i}}=k\ p+i'=p+i\;, \qquad k\in\nit\;, \quad 0\le i'\le p-1\;, \quad 1\le i\le\infty \;,\]
then the number of elements $n_e(p+r)$ of a complex completion $F_{\omega _{p+i}}$ \resp{$F_{\o\omega _{p+i}}$} at the infinite place $\omega _{p+i}$ \resp{$\o\omega _{p+i}$} is equal to:
\[ n_e(p+i) = m^{(p+i)}\cdot (p+i)\cdot N\cdot (n_{\rm nu})\;.\]

\item {\bbf The connection between the sets $\omega =\{\omega _1,\dots,\omega _r\}$ and $\o \omega =\{\o\omega _1,\dots,\o\omega _r\}$ of infinite places of $F$ and a set of finite places corresponding to finite extensions $K_p$ of $k_p$, for every prime $p$, can be obtained as follows:}

The field $K_p$, which is a finite extension of $k_p$ that may be $\QQ _p$, is a $p$-adic field.  Let $\Theta _{K_p}$ denote its ring of integers, $\wp_{K_p}$ the unique maximal ideal of $\Theta _{K_p}$, $k(v_{K_p})=k(\wp_{K_p})=\Theta _{K_p}/\wp_{K_p}$ its residue field, $\widetilde \omega _{K_p}$ a uniformizer in $\Theta _{K_p}$ and $v_{K_p}:K_p^*\to Z$ be the unique valuation.

The number of elements in $k(\wp_{K_p})$ is $q_p=p^{f_{v_{K_p}}}$ where $f_{v_{K_p}}=[k(v_{K_p}):\FF_p]$ is the residue degree over $k_p$.

On the other hand, let $E_p$ be an imaginary quadratic field in which $p$ splits.  Let  $F^+=F^+_R\cup F^+_L$ be the real field covering the complex field $F=F_R\cup F_L$ as described below.  If $F$ verifies $F=E_p\cdot F^+$, then $F$ is a CM field.

{\bbf A finite extension $K_p$ of $k_p$ can be identified with the completion of the real field $F^+$\/} in a place $v_{p+1}=v'_{p_1}$ above $p$. As, $p$ is decomposed in $E_p$ into two places $\wp$ and $\o \wp$, the places of $F$ dividing $p$ are divided into a set of left infinite complex places $\omega '_p=\{\dots,\omega'_{p_i},\dots,\omega '_{p_s}\}$ above $\wp$ and into a set of 
right infinite complex places $\o\omega '_p=\{\dots,\o\omega'_{p_i},\dots,\o\omega '_{p_s}\}$ above $\o \wp$ in such a way that:
\begin{align*}
\#\omega '_{p_i}&=\#\omega _{p+i}=p+i\;, \qquad 1\le i\le \infty \;,  \\
\rresp{
\#\o\omega '_{p_i}&=\# \o\omega _{p+i}=p+i},\end{align*}
where $(p+i)$ is in general given by $k_p+i'$.

By this way, {\bbf finite extensions $K_p$ can correspond to etale coverings of completions of $F$ in $\omega '_p$ and $\o\omega '_p$.}

\item Let $G^{(n)}(F_{\o\omega} \times F_\omega )$ be the bilinear algebraic semigroup with entries in the product, right by left, $F_{\o\omega} \times F_\omega $ of towers of packets of equivalent complex completions at the set 
$\o\omega \times\omega =\{\o\omega _1,\dots,\o\omega _r\} \times\{ \omega _1,\dots, \omega _r\}$ of biplaces.
 
 {\bbf $G^{(n)}(F_{\o\omega} \times F_\omega )$ is then composed of conjugacy class representatives $G^{(n)}(F_{\o\omega_{j,m_j}} \times F_{\omega_{j,m_j} })$, $1\le j\le r\le\infty $,\/} having multiplicities $m^{(j)}$, $1\le m_j\le m^{(j)}$, and {\bbf corresponding to the $r$ biplaces $(\o\omega \times\omega )$.}
 
Each conjugacy class representative $G^{(n)}(F_{\o\omega_{j,m_j} }\times F_{\omega_{j,m_j}} )$ is a $\GL_n(F_{\o\omega_{j,m_j}} \times F_{\omega_{j,m_j}} )$-subbisemimodule $\subset G^{(n)}(F_{\o\omega} \times F_\omega )$.

The decomposition of $ G^{(n)}(F_{\o\omega} \times F_\omega )$ into conjugacy classes can be realized by considering the cutting of the bilattice $\Lambda _{\o\omega }\times\Lambda _\omega $, referring to it, into subbilattices in $ G^{(n)}(F_{\o\omega} \times F_\omega )$ under the action of the product $T_R(n;r)\otimes T_L(n;r)$ of Hecke operators having as representation $\GL_n((\ZZ/N\ \ZZ)^2)$ \cite{Pie1}.

\item The algebraic representation of the bilinear algebraic semigroup of matrices $\GL_n(F_{\o\omega }\times F_\omega )$ into the 
$\GL_n(F_{\o\omega }\times F_\omega )$-bisemimodule 
$G^{(n)}(F_{\o\omega }\times F_\omega )$ corresponds to the algebraic morphism from
$\GL_n(F_{\o\omega }\times F_\omega )$ into
$\GL(G^{(n)}(F_{\o\omega }\times F_\omega ))$ which is the group of automorphisms of $G^{(n)}(F_{\o\omega }\times F_\omega )$.

\item Let $F_{\omega _\oplus}=\bigoplus_{j,m_j}F_{\omega _{j,m_j}}$
and $F_{\o\omega _\oplus}=\bigoplus_{j,m_j}F_{\o\omega _{j,m_j}}$
denote the sums of completions.

Then, $G^{(n)}(F_{\o\omega _\oplus}\times F_{\omega _\oplus})=
\bigoplus_{j,m_j} G^{(n)}(F_{\o\omega _{j,m_j}}\times F_{\omega _{j,m_j}})$ will correspond to the sums of conjugacy class representatives of the bilinear algebraic semigroup 
$G^{(n)}(F_{\o\omega }\times F_\omega )$.

In this respect, {\bbf $\GL(G^{(n)}(F_{\o\omega_\oplus }\times F_{\omega_\oplus} ))$ constitutes the $n$-dimensional equivalent of the product, right by left, $W^{ab}_{F_R} \times W^{ab}_{F_L}$ of the sums of the equivalence classes of the global Weil groups 
	$W^{ab}_{F_R}$ and $W^{ab}_{F_L}$ \cite{Pie1} and $
	G^{(n)}(F_{\o\omega _\oplus}\times F_{\omega _\oplus})$ becomes naturally its $n$-dimensional (irreducible) representation space according to}:
\[\Irr\Rep^{(n)}_{W_{F\RL}}(W^{ab}_{F_R}\times W^{ab}_{F_L})
=G^{(n)}(F_{\o\omega _\oplus}\times F_{\omega _\oplus})\]
if we have the injective morphism:
\[ W^{ab}_{F_R}\times W^{ab}_{F_L}\To 
\GL (G^{(n)}(F_{\o\omega _\oplus}\times F_{\omega _\oplus}))\;.\]
\Ei

\section{Cuspidal representations of the  bilinear algebraic semigroups}

\Bi
\item
The next step then consists in providing a (super)cuspidal representation of the bilinear algebraic semigroup $G^{(n)}(F_{\o\omega} \times F_\omega )$ or of $G^{(n)}(F_{\o\omega_\oplus }\times F_{\omega _\oplus})$.  The procedure can be summarized as follows:
\Bean
\item finding the cuspidal subrepresentation of each conjugacy class representative\linebreak $G^{(n)}(F_{\o\omega _{j,m_j}}\times F_{\omega_{j,m_j}} )$ of $G^{(n)}(F_{\o\omega} \times F_\omega )$.

\item showing that the sum of these cuspidal subrepresentations correspond to the searched cuspidal representation $\Irr\cusp(G^{(n)}(F_{\o\omega} \times F_\omega ))$ of the algebraic bilinear semigroup\linebreak $G^{(n)}(F_{\o\omega }\times F_\omega )$.
\Ee

\item Let 
\begin{align*}
\gamma ^T_{F_{\omega _{j,m_j}}}: \quad F_{\omega _{j,m_j}} &\To F^T_{\omega _{j,m_j}}\;,&& 1\le j\le r\le\infty \;, \\
\rresp{
\gamma ^T_{F_{\o\omega _{j,m_j}}}: \quad F_{\o\omega _{j,m_j}} &\To F^T                                                                                                                                                                     _{\o\omega _{j,m_j}}\;,}
\end{align*}
be the {\bbf toroidal isomorphism} mapping each \lr complex completion 
$F_{\omega _{j,m_j}}$ \resp{$F_{\o\omega _{j,m_j}}$} into its toroidal equivalent $F^T_{\omega _{j,m_j}}$ \resp{$F^T_{\o\omega _{j,m_j}}$} which is a complex one-dimensional semitorus $T^1_L[j,m_j]$ \resp{$T^1_R[j,m_j]$} localized in the upper \resp{lower} half space.

Then, {\bbf the morphisms:}
\begin{align*}
T_n(F^T_{\omega _{j,m_j}}): \quad F^T_{\omega _{j,m_j}} &\To T^{(n)}(F^T_{\omega _{j,m_j}})=T^n_L[j,m_j]\;,\\
\rresp{
T^t_n(F^T_{\o\omega _{j,m_j}}): \quad F^T_{\o\omega _{j,m_j}} &\To T^{(n)}(F^T_{\o\omega _{j,m_j}})=T^n_R[j,m_j]},\end{align*}
introduced in section 2, {\bbf sends the \lr toroidal complex completion 
$F^T_{\omega _{j,m_j}}$ \resp{$F^T_{\o\omega _{j,m_j}}$} into} the upper \resp{lower} semispace $T^{(n)}(F_{\omega _{j,m_j}})$ \resp{$T^{(n)}(F^T_{\o\omega _{j,m_j}})$} which is {\bbf a $n$-dimensional complex semitorus $T^n_L(j,m_j]$ \resp{$T^n_R[j,m_j]$}, corresponding to the conjugacy class representative $G^{(n)}(F^T_{\omega _{j,m_j}})$ \resp{$G^{(n)}(F_{\o\omega _{j,m_j}})$}.}

So, the composition of bimorphisms
\[ (T^t_n\circ T_n)\circ (\gamma ^T_{F_{\o\omega _{j,m_j}}}\times
 \gamma ^T_{F_{\omega _{j,m_j}}})
 :\quad
 F_{\o\omega _{j,m_j}} \times
 F_{\omega _{j,m_j}} \To
 G^{(n)}(F^T_{\o\omega _{j,m_j}}\times F^T_{\omega _{j,m_j}})\]
 \Bean
 \item is responsible for the generation of conjugacy class representatives\linebreak
 $G^{(n)}(F^T_{\o\omega _{j,m_j}}\times F^T_{\omega _{j,m_j}})$, $\forall\ j,m_j$ of the bilinear semigroup $G^{(n)}(F^T_{\omega }\times F^T_\omega )$ from the products, right by left, 
 $F_{\o\omega _{j,m_j}} \times F_{\omega _{j,m_j}}$
 of complex completions;
 
 \item allows to envisage the bimorphisms:
 \begin{multline*}
 \qquad\gamma ^T_{G^{(n)}(F_{\o\omega _{j,m_j}})}
 \times \gamma ^T_{G^{(n)}(F_{\omega _{j,m_j}})}: \quad
G^{(n)}(F_{\o\omega _{j,m_j}}
\times F_{\omega _{j,m_j}})\\
\To G^{(n)}(F^T_{\o\omega _{j,m_j}}
\times F^T_{\omega _{j,m_j}})\qquad\end{multline*}
which send each conjugacy class representative
$G^{(n)}(F_{\o\omega _{j,m_j}}
\times F_{\omega _{j,m_j}})$ of $G^{(n)}(F_{\o\omega }
\times F_{\omega })$ into its toroidal equivalent
$G^{(n)}(F^T_{\o\omega _{j,m_j}}
\times F^T_{\omega _{j,m_j}})$.
\Ee

\item Next, we consider the {\bbf \lr (semi)algebra $\widehat G^{(n)}(F^T_{\omega })$ \resp{$\widehat G^{(n)}(F^T_{\o\omega }$} }\cite{Pie1}, \cite{Pie2} {\bbf of continuous complex-valued measurable functions $\phi ^{(n)}_{G^T_L}(x_{g_L})$
\resp{$\phi ^{(n)}_{G^T_R}(x_{g_R})$} on $G^{(n)}(F^T_\omega )$
\resp{$G^{(n)}(F^T_{\o\omega} )$}} satisfying:
\begin{align*}
\int_{G^{(n)}(F^T_\omega )}\ \L|\phi ^{(n)}_{G^T_L}(x_{g_L})\R|\ dx_{g_L}&<\infty  \\
\rresp{
\int_{G^{(n)}(F^T_{\o\omega} )}\ \L|\phi ^{(n)}_{G^T_R}(x_{g_R})\R|\ dx_{g_R}&<\infty},
\end{align*}
with respect to a unique Haar measure on
$G^{(n)}(F^T_\omega )$
\resp{$G^{(n)}(F^T_{\o\omega} )$}: it is also noted $L^1_L(G^{(n)}(F^T_\omega ))$
\resp{$L^1_R(G^{(n)}(F^T_{\o\omega} )$}.

The bisemialgebra $\widehat G^{(n)}(F^T_{\o\omega}
\times F^T_{\omega} )$ of continuous complex valued measurable bifunctions
$\phi ^{(n)}_{G^T_R}(x_{g_R})\otimes \phi ^{(n)}_{G^T_L}(x_{g_L})$ on $G^{(n)}(F^T_{\o\omega}\times F^T_\omega )$ satisfying
\[ \int_{G^{(n)}(F^T_{\o\omega} \times F^T_\omega )}\ \L|\phi ^{(n)}_{G^T_R}(x_{g_R}) \otimes
\phi ^{(n)}_{G^T_L}(x_{g_L})\R|\ dx_{g_R}\ dx_{g_L}<\infty \;,\]
is noted $L^{1-1}\RL(G^{(n)}(F^T_{\o\omega }\times F^T_\omega ))$.

If the right functions $\phi ^{(n)}_{G^T_R}(x_{g_R})$ are projected on the left functions
$\phi ^{(n)}_{G^T_L}(x_{g_L})$, they become cofunctions and$L^{1-1}\RL(G^{(n)}(F^T_{\o\omega }\times F^T_\omega ))$ is transformed under this involution into the (bisemi)algebra
$L^2_{L\times L}(G^{(n)}(F^T_{\o\omega} \times F^T_\omega ))$ of square integrable functions.

\item Each \lr function $\phi ^{(n)}_{G^T_{j_L}}(x_{g_{j_L}})$
\resp{$\phi ^{(n)}_{G^T_{j_R}}(x_{g_{j_R}})$} of $L^1_L(G^{(n)}(F^T_\omega ))$ \resp{$L^1_R(G^{(n)}(F^T_{\o\omega }))$}
is defined on a conjugacy class representative
$G^{(n)}(F^T_{\omega _{j,m_j}})$
(resp.\linebreak {$G^{(n)}(F^T_{\o\omega _{j,m_j}})$}
of $G^{(n)}(F^T_\omega )$
\resp{$G^{(n)}(F^T_{\o\omega })$}
which is a $n$-dimensional complex \lr semitorus.

Thus,
$\phi ^{(n)}_{G^T_{j_L}}(x_{g_{j_L}})$
\resp{$\phi ^{(n)}_{G^T_{j_R}}(x_{g_{j_R}})$} is a \lr function
$\phi _L(T^n_L[j,m_j])$
\resp{$\phi _R(T^n_R[j,m_j])$} on the semitorus
\[ G^{(n)}(F^T_{\omega _{j,m_j}})\equiv T^n_L[j,m_j]
\rresp{G^{(n)}(F^T_{\o\omega _{j,m_j}})\equiv T^n_R[j,m_j]}.\]

\begin{propo} Each \lr function
 $\phi ^{(n)}_{G^T_{j_L}}(x_{g^T_{j_L}})$
\resp{$\phi ^{(n)}_{G^T_{j_R}}(x_{g^T_{j_R}})$} on the conjugacy class representative 
$G^{(n)}(F^T_{\omega _{j,m_j}})$
\resp{$G^{(n)}(F^T_{\o\omega _{j,m_j}})$} is a function on the \lr $n$-dimensional complex semitorus 
$T^n_L[j,m_j]$
\resp{$T^n_R[j,m_j]$} having the analytic development:
\begin{align*}
\phi _L(T^n_L[j,m_j])&=\lambda ^{\half}(n,j,m_j)\ e^{2\pi ijz}\\
\rresp{\phi _R(T^n_R[j,m_j])&=\lambda ^{\half}(n,j,m_j)\ e^{-2\pi ijz}},\end{align*}
where:
\Bi
\item $\vec z=\sum\limits^n_{d=1}z_d \ \vec e_d$ is a point of $G^{(n)}(F^T_{\omega _{j,m_j}})$;
\item $\lambda ^{\half}(n,j,m_j)\simeq j^n\ N^n\ (m^{(j)})^n$ can be considered as a Heche character.
\Ei
\end{propo}

\begin{proof}
\Be
\item As we are concerned with the $j$-th infinite complex place $\omega _j$ \resp{$\o\omega _j$}, we have to take into account the global Frobenius substitution given by the mapping:
\[ e^{2\pi iz}\To e^{2\pi ijz}
\rresp{e^{-2\pi iz}\To e^{-2\pi ijz}},\]
in such a way that
\begin{align*}
e^{e\pi ijz} & \simeq e^{2\pi ijz_1}\times \dots\times e^{2\pi ijz_n}\\
\rresp{e^{-2\pi ijz} & \simeq e^{-2\pi ijz_1}\times \dots\times e^{-2\pi ijz_n}}
\end{align*}
because the $n$-dimensional complex semitorus $T^n_L[j,m_j]$ \resp{$T^n_R[j,m_j]$} is diffeomorphic to the $n$-fold product:
\begin{align*}
T^n_L[j,m_j] &\simeq  T^1_L[j,m_j] \times \dots\times T^1_L[j,m_j]\\
\rresp{T^n_R[j,m_j] &\simeq  T^1_R[j,m_j] \times \dots\times T^1_R[j,m_j]};
\end{align*}

\item The scalar $\lambda (n,j,m_j)$ corresponds to the $(j,m_j)$-th coset representative $U_{j,m_{j_R}}\times U_{j,m_{j_L}}$ of the product
$T_R(n;r)\otimes T_L(n;r)$ of Hecke operators being represented by
$\GL_n((\ZZ/N\ \ZZ)^2)$.

More precisely, let $\{\lambda_d(n,j,m_j)\}^{2n}_{d=1}$ be the set of eigen(bi)values of $U_{j,m_{j_R}}\times U_{j,m_{j_L}}$ and let
$\lambda (n,j,m_j)=\prod\limits^{2n}_{d=1}\lambda_d(n,j,m_j)$ be the product of these eigenvalues.

According to \cite{Pie1}, we have that:
\[ \lambda (n,j,m_j) = \prod^{2n}_{d=1}\lambda _d(n,j,m_j)
=\det(\alpha _{n^2,j^2}\times D_{j^2,m_{j^2}})\simeq j^{2n}\cdot N^{2n}\cdot (m^{(j)})^{2n}\]
where:
\Bi
\item $D_{j^2,m_{j^2}}$ is the decomposition group of th $j$-th bisublattice with representative $m_j$;
\item $\alpha _{n^2nj^2}$ is the $j$-th split Cartan subgroup.
\Ei

It then appears that the square root $\lambda ^{\half}(n,j,m_j)$ of $\lambda (n,j,m_j)$ can be considered as a Hecke character having an inflation action on $e^{2\pi ijz}$.

\item Referring to the composition of $T^n_L[j,m_j]$ into a $n$-fold product of $T^1_L[j,m_j]$ and to the analytic development 
$\phi _L(T^n_L[j,m_j])=\lambda ^{\half}(n,j,m_j)\ e^{2\pi ijz}$ of the function $\phi _L$ on $T^n_L[j,m_j]$, it is clear that $\phi _L(T^1_L[,m_j])$ is a function on the $1$-dimensional complex semitorus $T^1_L[j,m_j]$.

Indeed, we have that:
\begin{align*}
\phi _L(T^1_L[j,m_j])
&= \lambda _{T^1_L}\ e^{2\pi ijz_d}  && z_d\in \CC\;, \\
&=S^1_{d_1}[j,m_j]\times S^1_{d_2}[j,m_j]\\
&= r_{S^1_{d_1}}\ e^{2\pi ijx_{d_1}}
\times r_{S^1_{d_2}}\ e^{2\pi ij(iy_{d_1})} && x_{d_1},y_{d_2}\in\rit\;, \\
&\simeq \lambda _{d_1}(n,j,m_j)
\times\lambda _{d_2}(n,j,m_j)\
e^{2\pi ijx_{d_1}}\cdot
e^{2\pi ij(iy_{d_2})}\end{align*}
where:
\Bi
\item $\lambda _{d_1}(n,j,m_j)
\cdot\lambda _{d_2}(n,j,m_j) \simeq r_{S^1_{d_1}}\cdot r_{S^1_{d_2}}$ is the product of the eigenvalues of $U_{j,m_{j_R}}\times U_{j,m_{j_L}}\subset T_R(1,r)\times T_L(1,r)$;
\item $r_{S^1_{d_1}}$ and $ r_{S^1_{d_2}}$ are radii of the circles
$S^1_{d_1}[j,m_j]$
and $S^1_{d_2}[j,m_j]$.
\Ei

Thus, the left $1D$-complex semitorus $T^1_L[j,m_j]$ is diffeomorphic to the product\linebreak
$S^1_{d_1}[j,m_j]\times S^1_{d_2}[j,m_j]$ of two circles localized in perpendicular planes with $\cos(2\pi ijy_{d_2})$ and
 $\sin(2\pi ijy_{d_2})$ of $e^{2\pi ij(iy_{d_2})}$ defined over $i\rit$.
\qedhere\Ee \end{proof}
 
\begin{propo}
As every \lr function $\phi _L(T^n_L[j,m_j])$
\resp{$\phi _R(T^n_R[j,m_j])$} constitutes the cuspidal representation
$\Pi _{j,m_j}(\GL_n(F_{\omega _{j,m_j}}))$
\resp{$\Pi _{j,m_j}(\GL_n(F_{\o\omega _{j,m_j}}))$}
of the $(j,m_j)$-th conjugacy class representative of the algebraic semigroup 
$\GL_n(F_\omega )$
\resp{$\GL_n(F_{\o\omega} )$}, the sum
$\bigoplus\limits_{j,m_j}(\phi _R(T^n_R[j,m_j])\otimes \phi _L(T^n_L[j,m_j]))$ of the cuspidal subrepresentations of all conjugacy class representatives of the bilinear algebraic semigroup $\GL_n(F_{\o\omega }\times F_\omega )$ is the searched cuspidal representation $\Pi(\GL_n(F_{\o\omega _\oplus}\times F_{\omega _\oplus}))$ according to:
\begin{align*}
\Pi(\GL_n(F_{\o\omega _\oplus}\times F_{\omega _\oplus}))
&= \bigoplus_{j,m_j}\Pi _{j,m_j}(\GL_n(F_{\o\omega _{j,m_j}}\times F_{\omega _{j,m_j}}))\\
&=\bigoplus_{j,m_j}(\phi _R(T^n_R[,m_j])\otimes \phi _L(T^n_L[j,m_j]))\;.\end{align*}
\end{propo}

\begin{proof}
If the sum $\bigoplus\limits^r_{j=1}$ tends to infinity, i.e. $r\to\infty$, then
$\bigoplus\limits_{j,m_j}\phi _L(T^n_L[j,m_j])$
\resp{$\bigoplus\limits_{j,m_j}\phi _R(T^n_R[j,m_j])$}
represents the Fourier development of a \lr cuspidal form over $\CC^n$.

And thus, $\Pi (\GL_n(F_{\o\omega _\oplus}\times F_{\omega _\oplus}))$ constitutes clearly the cuspidal representation of the bilinear algebraic semigroup $\GL_n(F_{\o\omega}\times F_{\omega })$, for $1\le j\le r\le\infty $.\end{proof}

\begin{propo}[Langlands global correspondence]
Let $\sigma _{j,m_j}(W_{F_{\o\omega _{j,m_j}}}\times W_{F_{\omega _{j,m_j}}})
= G^{(n)}(F_{\o\omega _{j,m_j}}\times F_{\omega _{j,m_j}})$ denote the $n$-dimensional representation subspace of the product, right by left,
$W_{F_{\o\omega _{j,m_j}}}\times W_{F_{\omega _{j,m_j}}}$ of the Weil subgroups restricted to
$F_{\o\omega _{j,m_j}}$ and  $F_{\omega _{j,m_j}}$ and given by
\[\sigma _{j,m_j}(W_{F_{\o\omega _{j,m_j}}}\times W_{F_{\omega _{j,m_j}}})
=\Irr\Rep^n (W_{F_{\o\omega _{j,m_j}}}\times W_{F_{\omega _{j,m_j}}})\]
as described in section 2.

Let 
\[\Pi _{j,m_j}(\GL_n(F_{\o\omega _{j,m_j}}\times {F_{\omega _{j,m_j}}})
=\Pi^\vee  _{j,m_j}(\GL_n(F_{\o\omega _{j,m_j}})\times  
\Pi   _{j,m_j}(\GL_n(F_{\omega _{j,m_j}})\]
be its cuspidal (sub)representation in such a way that
$\Pi^\vee  _{j,m_j}(\GL_n(F_{\o\omega _{j,m_j}}))$ be the contragradient cuspidal subrepresentation restricted to the Weil subgroup 
$W_{F_{\o\omega _{j,m_j}}}$.

Then, there exists bijective morphisms:
\begin{multline*}
\qquad T_{j,m_j}: \quad
\sigma _{j,m_j}(W_{F_{\o\omega _{j,m_j}}}\times W_{F_{\omega _{j,m_j}}})
\\ \To \Pi _{j,m_j}(\GL_n(F_{\o\omega _{j,m_j}}\times F_{\omega _{j,m_j}})\;, \qquad
1\le j\le r\le\infty \;,\qquad \end{multline*}
between the $n$-dimensional conjugacy class representatives of the products, right by left, of the Weil subgroups and the corresponding $n$-dimensional cuspidal class representatives,

leading to the {\bfseries\boldmath bijective morphism:}
\[ T:\quad (\sigma (W^{ab}_{F_R}\times W^{ab}_{F_L})\To
\Pi (\GL_n(F_{\o\omega }\times F_\omega ))\]
{\bfseries\boldmath between the sum $\sigma (W^{ab}_{F_R}\times W^{ab}_{F_L})$ of the $n$-dimensional conjugacy class representatives of the Weil subgroups given by the algebraic bilinear semigroup $G^{(n)}(F_{\o\omega _{\oplus}}\times F_{\omega _\oplus})$ and its cuspidal representation given by $\Pi (\GL_n(F_{\o\omega }\times F_\omega ))$.}
\end{propo}

\begin{proof}
\Be
\item The $n$-dimensional conjugacy class representative of the product, right by left, of the Weil subgroups 
$W_{F_{\o\omega _{j,m_j}}}\times W_{F_{\omega _{j,m_j}}}$ is given by:
\[ G^{(n)}({F_{\o\omega _{j,m_j}}}\times {F_{\omega _{j,m_j}}})
=\sigma _{j,m_j}(W_{F_{\o\omega _{j,m_j}}}\times W_{F_{\omega _{j,m_j}}})\;.\]
The toroidal compactification
\begin{align*}
T_{j,m_j}(G^{(n)}({F_{\o\omega _{j,m_j}}}\times {F_{\omega _{j,m_j}}}))
&\simeq \Pi _{j,m_j}(\GL_n({F_{\o\omega _{j,m_j}}}\times {F_{\omega _{j,m_j}}})\\
&=\lambda ^{\half}(n,j,m_j)\ e^{-2\pi ijz}\times
\lambda ^{\half}(n,j,m_j)\ e^{2\pi ijz}\end{align*}
of the conjugacy class representative
$\GL_n({F_{\o\omega _{j,m_j}}}\times {F_{\omega _{j,m_j}}})$ of the bilinear algebraic semigroup
$\GL_n( {F_{\o\omega }}\times {F_{\omega }})$ is in bijection with the corresponding cuspidal conjugacy class representative
$\Pi _{j,m_j}(\GL_n({F_{\o\omega _{j,m_j}}}\times {F_{\omega _{j,m_j}}}))$ given by
$\phi _R(T^n_R[j,m_j])\otimes \phi _L(T^n_L[j,m_j])$ as developed in proposition 3.1.

\item Then,  the sum
$\sigma (W^{ab}_{F_R}\times W^{ab}_{F_L})$ of the $n$-dimensional conjugacy class representatives of the Weil subgroups given by
\[
G^{(n)}(F_{\o\omega _{\oplus}}\times F_{\omega _\oplus})=
 \bigoplus_{j,m_j} G^{(n)}(F_{\o\omega _{j,m_j}}\times F_{\omega _{j,m_j}})\]
 is in one-to-one correspondence with the searched cuspidal representation\linebreak $\Pi (\GL_n(F_{\o\omega _{j,m_j}}\times F_{\omega _{j,m_j}})$ according to:
 \begin{align*}
 T(\sigma (W^{ab}_{F_R}\times W^{ab}_{F_L}))
 &\simeq
 \bigoplus_{j,m_j}\Pi _{j,m_j}(\GL_n(F_{\o\omega _{j,m_j}}\times F_{\omega _{j,m_j}})\\
 &=
\Pi (\GL_n(F_{\o\omega  }\times F_{\omega  }).\qedhere
\end{align*}
\Ee
\end{proof}

\Ei
\section{The Langlands functoriality in this new bilinear mathematical framework}

\Bi
\item The Langlands global correspondence(s) having been stated in the irreducible complex case, it is now time to ask {\bbf in what extent the cuspidal representation $\Pi ^{(2n)}(\GL_{2n}(F_{\o\omega }\times F_\omega ))$ of the algebraic bilinear semigroup $\GL_{2n}(F_{\o\omega }\times F_\omega )$ can be reached from the knowledge of the cuspidal representation $\Pi ^{(2)}(\GL_{2}(F_{\o\omega }\times F_\omega ))$ of the algebraic bilinear semigroup $\GL_{2}(F_{\o\omega }\times F_\omega )$.}  

As recalled in section 1, the Langlands functoriality conjecture consists in proving that:
\[
\Pi ^{(2n)}
= \sym^{(n)}(\Pi ^{(2n)})
= \mathop\otimes\limits_{v_p}\sym^{(n)}(\Pi ^{(2)}_{v_p})\]
where $\Pi ^{(2)}_{v_p}$ denotes the cuspidal subrepresentation of the considered algebraic group at the primary place $v_p$.

{\bbf This case of functoriality, transposed in the considered bilinear framework, will be considered  in the following, but also a case of functoriality associated with the cuspidal reducible representation of $\GL_{2n}(F_{\o\omega }\times F_\omega )$, which constitutes the content of our main proposition 4.2.}

\begin{propo}[\bbf Functoriality $\sym^{(n)}$ in a bilinear framework]
Let 
\[ \sym^{(n)}: \quad \GL_2(F_{\o\omega }\times F_\omega )
\To \GL_{2n}(F_{\o\omega }\times F_\omega )\]
or equivalently:
\[ \sym^{(n)}: \quad G^{(2)}(F_{\o\omega }\times F_\omega )
\To G^{(2n)}(F_{\o\omega }\times F_\omega )\]
be the $n$-th symmetric power representation of the algebraic bilinear semigroup\linebreak $\GL_2(F_{\o\omega }\times F_\omega )
$.

Then, {\bbf the cuspidal representation 
$\Pi ^{(2n)}(\GL_{2n}(F_{\o\omega }\times F_\omega ))$ of the bilinear algebraic semigroup $\GL_{2n}(F_{\o\omega }\times F_\omega )$ can be reached functorially from the cuspidal representation
$\Pi ^{(2)}(\GL_{2}(F_{\o\omega }\times F_\omega ))$ of the algebraic bilinear semigroup 
$\GL_{2}(F_{\o\omega }\times F_\omega )$ throughout the injective morphism}:
\[\Pi ^{\text{\rm cusp}}_{2\to 2n}: \quad
\Pi ^{(2)}(\GL_{2}(F_{\o\omega }\times F_\omega ))
\To \Pi ^{(2n)}(\GL_{2n}(F_{\o\omega }\times F_\omega ))\]
which corresponds to the morphism $\sym^{(n)}$ on the cuspidal representation
$\Pi ^{(2)}(\GL_{2}(F_{\o\omega }\times F_\omega ))$
of $\GL_{2}(F_{\o\omega }\times F_\omega )$.
\end{propo}

\begin{proof}
\Be
\item Referring to proposition 3.3, we have that the cuspidal representation\linebreak $\Pi ^{(2)}(\GL_{2}(F_{\o\omega }\times F_\omega ))$ of the bilinear algebraic semigroup
$\GL_{2}(F_{\o\omega }\times F_\omega )$ can be developed according to:
\begin{multline*}
\begin{aligned}
\Pi ^{(2)}(\GL_{2}(F_{\o\omega }\times F_\omega ))
&= \bigoplus_{j,m_j}
\Pi ^{(2)}_{j,m_j}(\GL_{2}(F_{\o\omega_{j,m_j} }\times F_{\omega_{j,m_j}} ))\\
&= \bigoplus_{j,m_j} (\phi _R(T^2_R[j,m_j])\otimes \phi _L(T^2_L[j,m_j]))
\\
&=\bigoplus_{j,m_j}(\lambda ^{\half}(2,j,m_j)\ e^{-2\pi ijy_2}
\otimes\lambda ^{\half}(2,j,m_j)\ e^{2\pi ijy_2} \;,
\end{aligned}\\
  y_2\in G^{(2)}(F_\omega )\;,
\end{multline*}
with respect to the cuspidal subrepresentations of the conjugacy class representatives $\GL_{2}(F_{\o\omega_{j,m_j} }\times F_{\omega _{j,m_j}}))$ of $\GL_{2}(F_{\o\omega }\times F_\omega )$.

Remark that the complex dimensions are here envisaged in real notations.

\item Then, every cuspidal conjugacy class representative
$\Pi ^{(2n)}_{j,m_j}(\GL_{2n}(F_{\o\omega_{j,m_j} }\times F_{\omega_{j,m_j}} ))$ of the bilinear algebraic semigroup
$\GL_{2n}(F_{\o\omega  }\times F_{\omega } )$ can be obtained from the corresponding cuspidal conjugacy class representative 
$\Pi ^{(2)}_{j,m_j}(\GL_{2}(F_{\o\omega_{j,m_j} }\times F_{\omega_{j,m_j}} ))$  of
$\GL_{2}(F_{\o\omega}\times F_{\omega} )$ by means of the injective morphism:
\begin{multline*}
\Pi ^{\text{\rm cusp}}_{2\to 2n}(j,m_j): \quad
\lambda ^{\half}(2,j,m_j)\ e^{-2\pi ijy_2}\otimes
\lambda ^{\half}(2,j,m_j)\ e^{2\pi ijy_2}\\
\qquad \qquad\To \lambda ^{\half}(2n,j,m_j)\ e^{-2\pi ijy}\otimes
 \lambda ^{\half}(2n,j,m_j)\ e^{2\pi ijy}\;, \\
 y_2\in G^{(2)}(F_\omega )\;, 
 y\in G^{(2n)}(F_\omega )\;, \end{multline*}
 sending:
 \Bi
 \item $\lambda ^{\half}(2,j,m_j)$ into $\lambda ^{\half}(2n,j,m_j)$
\item $y_2$ into $y$.
\Ei
This is possible if proposition 3.1 is taken into account.

\item And, the searched cuspidal representation
$\Pi ^{(2n)}(\GL_{2n}(F_{\o\omega  }\times F_{\omega } ))$ 
results from the sum of the injective morphisms:
\begin{multline*}
\bigoplus_{j,m_j}\L[ 
\Pi ^{\text{\rm cusp}}_{2\to 2n}(j,m_j): \quad
\Pi ^{(2)}_{j,m_j}(\GL_2( F_{\o\omega _{j,m_j}}\times F_{\omega _{j,m_j}}))\R.\\
\To\L.
\Pi ^{(2n)}_{j,m_j}(\GL_{2n}( F_{\o\omega _{j,m_j}}\times F_{\omega _{j,m_j}}))\R]
\end{multline*}
in such a way that
\[
\Pi ^{(2n)} (\GL_{2n}( F_{\o\omega  }\times F_{\omega  }))
=\bigoplus_{j,m_j}
\Pi ^{(2n)}_{j,m_j}(\GL_{2n}( F_{\o\omega _{j,m_j}}\times F_{\omega _{j,m_j}}))\;.\qedhere\]
\Ee\end{proof}

\item This treatment of functoriality is rather trivial in the considered bilinear framework.

{\bbf A more interesting way of envisaging functoriality is to take into account the fact that the algebraic bilinear semigroup $\GL_2(F_{\o\omega }\times F_\omega )$ is a bisemigroup submitted to the cross binary operation ${\u\times}$ which allows to reduce
the problem of Langlands functoriality to the reducibility of representations of groups.}

Let $2n=2_1+
+2_2+\dots+2_\ell +\dots+2_n$ be a partition of the integer $2n$ and let
\[ \GL_{2_1}(F_{\o\omega }\times F_\omega )
\times \GL_{2_2}(F_{\o\omega }\times F_\omega )
\times\dots\times
\GL_{2_\ell }(F_{\o\omega }\times F_\omega )
\times\dots\times
\GL_{2_n }(F_{\o\omega }\times F_\omega )\]
be the $n$-th symmetric power of $\GL_{2}(F_{\o\omega }\times F_\omega )$ according to this partition.

Referring to section 2 and to \cite{Pie2}, it appears that the product ``$\times$'' between two bilinear algebraic semigroups is the cross binary operation ``$\u\times$'' which enables to develop this product according to:
\begin{multline*}
 \GL_{2_1}(F_{\o\omega }\times F_\omega )
\u\times\dots\u\times
GL_{2_\ell }(F_{\o\omega }\times F_\omega )
\u\times\dots\u\times
\GL_{2_n }(F_{\o\omega }\times F_\omega )\\
= (\GL_{2_1}(F_{\o\omega })
\oplus\dots\oplus
\GL_{2_\ell }(F_{\o\omega })
\oplus\dots\oplus
\GL_{2_n }(F_{\o\omega }))\\
\times
(\GL_{2_1}(F_{\omega })
\oplus\dots\oplus
\GL_{2_\ell }(F_{\omega })
\oplus\dots\oplus
\GL_{2_n }(F_{\omega }))
\end{multline*}
and to state the main proposition.

\begin{Mpropo}
The cuspidal representation $\Pi ^{(2n)} (\GL_{(2n)}( F_{\o\omega  }\times F_{\omega  }))$ of the bilinear algebraic semigroup
$\GL_{2n}( F_{\o\omega  }\times F_{\omega  })$ is (non orthogonally) completely reducible if it decomposes diagonally according to the direct sum $\bigoplus\limits_{\ell=1}^n\Pi ^{(2_\ell)}
(\GL_{2_\ell }( F_{\o\omega  }\times F_{\omega  }))$ of irreducible cuspidal representations of the algebraic bilinear semigroups
$\GL_{2_\ell }( F_{\o\omega  }\times F_{\omega  })$ and off-diagonally according to the direct sum
$\bigoplus\limits_{k\neq\ell}(
\Pi ^{(2_k)} (\GL_{2_k }( F_{\o\omega  }))
\otimes\Pi ^{(2_\ell )} (\GL_{2_\ell  }( F_{\omega  })))$ of the (tensor) products of irreducible cuspidal representations of cross algebraic linear semigroups
$\GL_{2_k}(F_{\o\omega })
\times\GL_{2_\ell}(F_{\omega })\equiv T^t_{2_k}(F_{\o\omega} )\times T_{2_\ell}(F_\omega )$, $\forall\ k\neq \ell$, $1\le k,\ell\le n$.

This reducible cuspidal representation
\[ 
\Pi ^{(2n)} (\GL_{2n}( F_{\o\omega  }\times F_{\omega  }))
= \bigoplus_{\ell=1}^n \Pi ^{(2_\ell)}(\GL_{2_\ell}( F_{\o\omega  }\times F_{\omega  }))
 \bigoplus_{k\neq\ell=1}^n (\Pi ^{(2_k)}(\GL_{2_k}( F_{\o\omega  }))
 \otimes \Pi ^{(2_\ell)}(\GL_{2_\ell}( F_{\omega  })))\]
 of $\GL_{2n}( F_{\o\omega  }\times F_{\omega  })$ then corresponds to the Langlands functoriality:
 \[ 
\Pi ^{\text{\rm cusp}}_{2\to 2n}: \quad
\Pi ^{(2)}(\GL_2( F_{\o\omega }\times F_{\omega }))
\To
\Pi ^{(2n)} (\GL_{2n}( F_{\o\omega  }\times F_{\omega}  ))\;
\]
\end{Mpropo}

\begin{proof}
\Be
\item As it was noticed above, the cross binary operation ``$\u\times$'' allows to develop the $n$-th symmetric power of $\GL_2(F_{\o\omega }\times F_\omega )$ according to:
\begin{multline*}
 \GL_{2_1}(F_{\o\omega }\times F_\omega )
\times\dots\times
\GL_{2_\ell }(F_{\o\omega }\times F_\omega )
\times\dots\times
\GL_{2_n }(F_{\o\omega }\times F_\omega )\\
\begin{aligned}
&= \GL_{2_1}(F_{\o\omega }\times F_\omega )
\u\times\dots\u\times
 \GL_{2_\ell }(F_{\o\omega }\times F_\omega )
\u\times\dots\u\times
 \GL_{2_n }(F_{\o\omega }\times F_\omega )\\
&= \L( \bigoplus_{\ell=1}^n\GL_{2_\ell}(F_{\o\omega }\R)
\times\L( \bigoplus_{\ell=1}^n\GL_{2_\ell}(F_{\omega }\R)\\
&= \bigoplus_{\ell=1}^n\GL_{2_\ell}(F_{\o\omega }\times F_\omega )
\bigoplus_{k\neq \ell=1}^n(T^t_{2_k} (F_{\o\omega })\times T_{2_\ell}(F_\omega ))\;.
\end{aligned}\end{multline*}
The sum of the cuspidal representations of these algebraic bilinear semigroups\linebreak $\GL_{2_\ell}(F_{\o\omega }\times F_\omega )$ and $(T^t_{2_k} (F_{\o\omega })\times T_{2_\ell}(F_\omega ))$ is the reducible cuspidal representation of the algebraic bilinear semigroup $\GL_{2n}(F_{\o\omega }\times F_\omega )$ according to:
\begin{multline*} 
\Pi ^{(2n)} (\GL_{2n}( F_{\o\omega  }\times F_{\omega} ))\\
= \bigoplus_{\ell=1}^n\Pi ^{(2_\ell)}(\GL_{2_\ell}(F_{\o\omega }\times F_\omega ))
 \bigoplus_{k\neq\ell=1}^n
\L(\Pi ^{(2_k)}(\GL_{2_k}(F_{\o\omega }))\times
\Pi ^{(2_\ell)}(\GL_{2_\ell}(F_{\omega }))\R)\;.\end{multline*} 

\item This decomposition of the cuspidal representation of 
$\GL_{2n}( F_{\o\omega  }\times F_{\omega} )$ then corresponds clearly to the Langlands functoriality statement:
\[ 
\Pi ^{\text{\rm cusp}}_{2\to 2n}: \quad
\Pi ^{(2)}(\GL_2( F_{\o\omega }\times F_{\omega }))
\To
\Pi ^{(2n)} (\GL_{2n}( F_{\o\omega  }\times F_{\omega}  ))\;\qedhere
\]
\Ee \end{proof}

\begin{coro}
The cuspidal representation $\Pi ^{(2n)} (\GL_{2n}( F_{\o\omega  }\times F_{\omega} ))$ of the bilinear algebraic semigroup
$\GL_{2n}( F_{\o\omega  }\times F_{\omega} )$ is orthogonally completely reducible if it is decomposed only diagonally according to the direct sum of the irreducible cuspidal representations of the algebraic bilinear semigroups
$\GL_{2_\ell}( F_{\o\omega  }\times F_{\omega} )$, $1\le\ell\le n$, as follows:
\[ \Pi ^{(2n)} (\GL_{2n}( F_{\o\omega  }\times F_{\omega} ))
= \bigoplus_{\ell=1}^n \Pi ^{(2_\ell)} (\GL_{2_\ell}( F_{\o\omega  }\times F_{\omega} ))\;.\]
\end{coro}

\begin{proof} This cuspidal representation of
$\GL_{2n}( F_{\o\omega  }\times F_{\omega} )$ is orthogonally completely reducible with respect to the non orthogonally completely reducible cuspidal representation considered in proposition 4.2 in the sense that the off-diagonal cuspidal representations\linebreak
$\Pi ^{(2_k)} (\GL_{2_k}( F_{\o\omega  }))
\otimes \Pi ^{(2_\ell)} (\GL_{2_k}( F_{\omega  }))$ are not taken into account.\end{proof}

\Ei

       }
\end{document}